\documentclass[12pt]{article}
\usepackage{amssymb}
\usepackage{enumerate}
\usepackage{graphicx}

\def\R{{\rm I\! R}}

\def\C{\mbox{l\hspace{-.47em}C}}

\def\psx{{\partial p \over \partial x}}
\def\qsx{{\partial q \over \partial x}}
\def\psy{{\partial p \over \partial y}}
\def\qsy{{\partial q \over \partial y}}
\def\psu{{\partial p^* \over \partial u}}
\def\qsu{{\partial q^* \over \partial u}}
\def\psv{{\partial p^* \over \partial v}}
\def\qsv{{\partial q^* \over \partial v}}

\newtheorem{theorem}{Theorem}
\newtheorem{corollary}{Corollary}

\newtheorem{lemma}{Lemma}

\title{A note on the divergence-free Jacobian Conjecture in $\R^2$
\thanks{Partially supported by the PRIN group \lq\lq Dinamica anolonoma,
hamiltoniana, delle popolazioni e dei sistemi piani\rq\rq and by the
GNAMPA group  \lq\lq Analisi qualitativa e comportamento asintotico
di equazioni differenziali ordinarie e di equazioni alle differenze
\rq\rq. 
2000 AMS {\it Mathematics Subject Classification }: 14R15, 34C08.
\ \ \ \ \ \ \ \ \ \ \ \ \ \ \ \ \ \ \ \ \ \ \ \ \ \ \ \ \ \ \ \ \ \ \ \ \ \ \ \ \ \ \ \ \ \ \ \ \ \ \ \ \ \ \ \ \ \ \ \ \ \ \ \ \ \ \ \ \ \ \ \ \ \ \ \ \ \ \ \ \ 
{\it Key words and phrases. }  Jacobian Conjecture, divergence-free.}
}

\author{M. Sabatini}
\begin{document}
\maketitle
\begin{abstract} We give a shorter proof to a recent result by Neuberger \cite{N}, in the real case.  Our result is essentially an application of the global asymptotic stability Jacobian Conjecture. We also extend some of the results presented in \cite{N}. \end{abstract}

\section{Introduction}

The classical Jacobian Conjecture was formulated in \cite{K} as a problem about the global invertibility of polynomial maps $\Phi:\C^n \mapsto \rightarrow \C^n$. Keller asked whether a polynomial map with constant non-zero Jacobian determinant is globally invertible, and its inverse is itself a polynomial map. The problem was widely studied in subsequent decades, producing several partial results and even some faulty proofs. In \cite{BCW} one finds a historical overview of research about the Jacobian Conjecture and a  rich survey of results published up to 1982. The paper \cite{dBvdE} contains a more recent list of results and some equivalent formulations of the problem in arbitrary dimension.  Among general results concerning such a problem, it is known that it is equivalent to prove or disprove the statement in any field of zero characteristic, that it is sufficient to prove $\Phi$'s injectivity in order to get its surjectivity \cite{BCW}, and that $\Phi$'s global invertibility implies that $\Phi^{-1}$ is a polynomial map. 
The most studied special case is the bidimensional one, $\Phi (x,y) = (P(x,y) ,Q(x,y) )$, where the statement was proved under the hypothesis that either P's or Q's degree is 4, or prime, or both degrees are $\leq$ 100 (see \cite{BCW} for a more comprehensive list of results). A recent result, which is the object of this paper, proves the global invertibility of jacobian maps of the form $\Phi (x,y) = (x+p(x,y),y+q(x,y))$, $p(x,y)$ and $q(x,y)$ without terms of degree 1, under the additional assumptions that $\psx + \qsy =0$ and $\psx \qsy - \psy \qsx  =0$. 
In higher dimensions, a striking result states that, in order to prove the $n$-dimensional Jacobian conjecture, it is sufficient to prove it for maps of the form $\Phi = L + C$, $L$ linear, $C$ cubic,  \cite{BCW}, or even for maps of the form $\Phi(X) = X + (AX)^3$, where A is a nilpotent matrix \cite{D}. 

A different question, arising in differential equations from the study of a critical point's global stability, is also known as a Jacobian Conjecture.  It is concerned with the global asymptotic stability (g. a. s.) of a critical point of a vector field whose jacobian eigenvalues have negative real part at every point of the space \cite{MY}. In \cite{O} it was showed that under such hypotheses, it is equivalent to prove the global asymptotic stability of a critical point or the global injectivity of the vector field. Such a result gave a new direction to the research about the g. a. s. Jacobian Conjecture. Thanks also to such a new approach, such a question was positively settled in dimension 2 in \cite{F}, \cite{Gl}, \cite{Gu}. In higher dimensions it is known to be false \cite{CEGHM}, unless some additional hypotheses hold. 

In this paper we give a shorter proof to Neuberger's result \cite{N} in the real case, showing that it is actually a consequence of the bidimensional g. a. s.  Jacobian Conjecture. Actually, we prove something more, since we do not make assumptions on the terms in $\psx \qsy - \psy \qsx $. Actually, for jacobian maps it is sufficient to require that $\psx + \qsy \geq 0$. Then we look for algebraic-like conditions which imply such a property, involving the degree and the order of the real polynomials $P$ and $Q$, or the degrees of the monomials contained in $P$ and $Q$. We also extend some of  the corollaries proved in \cite{N}, weakening some symmetry conditions. 

\section{Results}

Througout this paper we only consider polynomials with real coefficients. Given a polynomial $P$, we write $d(P)$ for its degree, $o(P)$ for its order.  We say that a polynomial is {\it even} if it is the sum of even-degree monomials,  {\it odd} if it is the sum of odd-degree monomials. 
Similalrly, we say that a polynomial is {\it x-even} if it contains only terms with even powers of $x$, {\it x-odd} if it contains only terms with odd powers of $x$.

We say that a non-negative integer is a {\it gap} of $P$ if it is the difference of the degrees of two distinct monomials in $P$. We denote by $G(P)$ the gap-set of $P$. As an example, the polynomial $P(x,y) = x^3 + y^3 + x^2y^2 + y^7$ has gap-set $G(P) = \{0, 1, 3, 4 \}$. If $P$ has exactly one monomial, then we say that it has empty gap-set. 

We say that the couple of polynomials $(P,Q)$ satisfies the {\it gap condition} if for every monomial $M$ in $P$, one has $d(M) -1 \not\in G(Q)$. The gap condition is not symmetric, as shown by the couple $(P,Q) = (x+y^2,x^6+y^2)$. In such a case one has $G(P) = \{ 1 \} $, $G(Q) = \{ 4 \} $, so that $(P,Q)$ satisfies the gap condition, but $(Q,P)$ does not.

We say that  $(P,Q)$ satisfies the {\it symmetric gap condition} if both  $(P,Q)$ and $(Q,P)$ satisfy the gap condition. 

Let $\Phi : \R^2 \rightarrow \R^2$, $\Phi(x,y)\equiv (P(x,y),Q(x,y))$ be a real polynomial map. Let $ J_\Phi$ be its jacobian matrix. We say that $\Phi $ is a {\it jacobian map} if its jacobian determinant $\det J_\Phi$ is a non-zero constant.  We first consider a straightforward consequence of the g. a. s. Jacobian Conjecture.

\begin{lemma}\label{lemma}
 Let $\Phi^*:\R^2 \rightarrow \R^2$, $\Phi^* (u,v) =  (u+p^*(u,v),v+q^*(u,v))$ be a jacobian polynomial map, with $o(p^*) > 1$, $o(q^*) > 1$. If  $\psu + \qsv \geq 0$, then $\Phi^*$ is injective.
\end{lemma}
{\it Proof.}
Since $\det J_{\Phi^*}$ is constant, its value can be evaluated  at the origin, hence $\det J_{\Phi^*} = 1$.
Let us consider the planar differential system associated to the map $-\Phi^*$,
$$
\dot u = - u - p^*(u,v), \qquad \dot v = - v - q^*(u,v).
$$
Its jacobian matrix is $-J_{\Phi^*}$. Since $\psu + \qsv \geq 0$, at every point of the plane $-J_{\Phi^*}$ has trace $ \leq -2 < 0$, and  determinant $\det (-J_{\Phi^*} )= \det J_{\Phi^*} =1 > 0$, hence its eigenvalues have negative real part at every point of $\R^2$. Since the g. a. s. Jacbian Conjecture holds, $-\Phi^*$ is injective, hence $\Phi^*$ is injective, too.  
\hfill$\clubsuit$\bigskip

Neuberger's result for real maps is contained in lemma \ref{lemma}, since in \cite{N} only maps with $\psu + \qsv = 0$ are considered. 

In relation to the  classical Jacobian Conjecture, algebraic-like hypotheses are usually considered. In fact, even if checking whether  $\psu + \qsv \geq 0$ in some cases can be done,  statements related to the map's degree or order are desirable. 

Condition $ii)$ of next theorem has been added only because we do not assume the linear part of $\Phi$ to be the identity, but the argument is the same as in \cite{N}. The other conditions, as well as those ones in theorem \ref{notmain}, are new.

\begin{theorem}\label{main}
Let $\Phi: \R^2 \rightarrow \R^2$ be a jacobian map of the type $\Phi(x,y) = (ax+by+p(x,y),cx+dy+q(x,y))$, $a,b,c,d \in \R$, $o(p) > 1$, $o(q) > 1$. If one of the following holds, \\
\indent i)   $\max \{d(p),d(q)\} <  o(p) + o(q) -1$,  \\
\indent  ii)  both $p(x,y)$ and $q(x,y)$ are even polynomials,  \\  
\indent iii) $p$ is odd, $q$ is even and $(p,q)$ satisfies the  gap condition,  \\
\indent iv) $(p,q)$ satisfies the symmetric gap condition,  \\
\noindent  then  $\Phi$ is globally invertible.
\end{theorem} 
{\it Proof.}
Without loss of generality, we may assume $ad - bc > 0$. Let $A$ be the linear map associated to the matrix
$$
\left(  
\begin{array}{rl}
 a & b  \\ c & d
\end{array} 
\right) .
$$
Let $ A^{-1}$ be its inverse. Let us set $\Phi^*(u,v) = \Phi(A^{-1}(u,v))$. The linear part of 
$\Phi^*$ is just the composition of $A^{-1}$ and $A$, hence it is the identity. Then, one has $\Phi^*(u,v) = (u + p^*(u,v),v + q^*(u,v))$, with $\det J_{\Phi^*} > 0$. A linear change of variables does not change a polynomial's order, its degree and the property of being even or odd, as above  defined. Hence $o(p^*) > 1$, $o(q^*) > 1$, and conditions $i)$, \dots, $iv)$ hold for $p^*$ and $q^*$ as well. Moreover, one has $\det J_{\Phi^*} > 0 $. Without loss of generality we may assume $\det J_{\Phi^*} = 1 $.

Computing the jacobian determinant of $\Phi^*$ gives
$$
\det J_{\Phi^*} = 1 + \left( \psu + \qsv \right) +\left(  \psu \qsv - \psv \qsu \right) .
$$
Let us set
$$
T^* =  \psu + \qsv , \qquad D^*  =  \psu \qsv - \psv \qsu   .
$$
Since $o(T^* +  D^*) > 0$, one has $T^* +  D^* \equiv  0$. 

In order to prove $i)$, consider that  the highest degree monomial in $T^*$ has degree $\leq \max \{d(p),d(q)\}  -1$, while the lowest degree monomial in $D^*$ has degree $\geq o(p) + o(q) -2$. If $\max \{d(p),d(q)\} <  o(p) + o(q) -1$, then $T^*$ and $D^*$ have no monomials of the same degree, hence, from $T^* +  D^* \equiv  0$, one gets both $T^* \equiv  0$ and $ D^* \equiv  0$. Applying the lemma \ref{lemma} one proves that $-\Phi^*$ is injective, hence $\Phi^*$ and $\Phi$ are injective, too. 

Now, in order to prove $ii)$, consider that if both $p(x,y)$ and $q(x,y)$ are even polynomials, then $p^*(x,y)$ and $q^*(x,y)$ are even, $T^*$ is odd and $D^*$ is even. From $T^* +  D^* \equiv  0$, one gets again $T^* \equiv  0$ and $ D^* \equiv  0$, since monomials of $T^*$ do not cancel with monomials of $D^*$. Then one can proceed as in point $i)$ for the  injectivity of $\Phi$.

Under the hypotheses of $iii)$,   $\psu$ is even, $\qsu$ is odd, $D^*$ is odd. Assume by absurd that there exists a positive integer $h$ such that both $\qsu$  and $D^*$ have a monomial of  degree $h$. Then $q$ has a monomial $M$ such that $d(M) -1 = h$. Also, there exist monomials $K$ in $p$ and $L$ in $q$ such that $d(K) + d(L) - 2 = h$. Hence $d(M) - d(L) =  d(K) - 1$, contradicting the gap condition. This proves that $T^*$ and $D^*$ have no monomials of the same degree, so that  $T^* \equiv  0$ and $ D^* \equiv  0$. Then the above argument applies.

Finally, if $iv)$ holds, assume by absurd that there exists a positive integer $h$ such that both $T^*$  and $D^*$ have a monomial of  degree $h$. Then, either $p$ or $q$ has a monomial $M$ such that $d(M) -1 = h$. If $M$ is in $q$, we may repeat the argument of point $iii)$.  If $M$ is in $p$, we may repeat the argument of point $iii)$, exchanging the roles of $p$ and $q$. 
\hfill$\clubsuit$\bigskip

In order to show that we are considering non-empty hypotheses, we give some examples of jacobian mappings satisfying the above conditions. For condition $i)$ we may choose Meisters' maps,
$$
\Phi(x,y) = (ax+by + \mu(\alpha a + \beta b) (\alpha y - \beta x)^2, cx + dy + 
\mu(\alpha c + \beta d)(\alpha y - \beta x)^2),
$$
with $\mu \neq 0$, $(\alpha , \beta )\neq (0,0)$, $ad - bc \neq 0$. For condition $ii)$ we may consider
 $$
\Phi(x,y) = (x+y + x^5 + x^6, y + x^5 + x^6).
$$
An example of map satisfying both conditions $iii)$ and $iv)$ is 
$$
\Phi(x,y) = (x+y^3,y).
$$

In \cite{BCW} it was proved that in order to prove the Jacobian Conjecture it is sufficient to prove it for cubic-linear jacobian maps. It may be interesting to show that in $\R^2$ every jacobian map of the type linear + homogeneous is invertible. 

\begin{corollary}\label{hom}
Let $\Phi: \R^2 \rightarrow \R^2$ be a jacobian map of the type $\Phi(x,y) = (ax+by+p_n(x,y),cx+dy+q_n(x,y))$, $a,b,c,d \in \R$, with $p_n$ and $q_n$ homogeneous polynomials of the same degree $n > 1$. Then  $\Phi$ is globally invertible.
\end{corollary} 
{\it Proof.}
One has  $ o(p) = o(q) = d(p) = d(q) = n$, hence condition $i)$ is satisfied: $\max \{d(p),d(q)\} = n < 2n - 1 = o(p) + o(q) -1$.
\hfill$\clubsuit$\bigskip

An example of planar  linear + cubic jacobian map is the following map,
$$
P(x,y) = 2x-y+x^3+x^2y+\frac {xy^2}3+\frac {y^3}{27},
\quad Q := 3x-3y+ \frac {12 x^3 }{ 5}+ \frac {12 x^2y}{ 5} + \frac { 4xy^2}{5 } + \frac { 4 y^3}{45 }.
$$

A result slightly different form theorem \ref{main} can be proved by assuming different symmetries on $p$ and $q$. In next statement we assume $\Phi$ to be of the form $(x+p(x,y), y+q(x,y))$, as in \cite{N}, since a linear change of variables in general does not preserve the requested symmetry property. 

\begin{theorem}\label{notmain}
Let $\Phi: \R^2 \rightarrow \R^2$ be a jacobian map of the type $\Phi(x,y) = (x+p(x,y), y+q(x,y))$, $o(p) > 1$, $o(q) > 1$. If one of the following holds, \\  
\indent i)  $p$ is $x$-even, $q$ is $x$-odd, \\  
\indent  ii) $p$ is $y$-odd, $q$ is $y$-even, \\
\noindent  then  $\Phi$ is globally invertible.
\end{theorem} 
{\it Proof.}
Working as in theorem \ref{main}, one has 
$$
\det J_{\Phi} = 1 + \left( \psx + \qsy \right) +\left(  \psx \qsy - \psy \qsx \right) .
$$
Let us set
$$
T =  \psx + \qsy  , \qquad D  = \psx \qsy - \psy \qsx   .
$$
As in theorem \ref{main}, since $o(T +  D) > 0$, one has $T +  D \equiv  0$. 

If $i)$ holds, then $T$ is $x$-odd, $D$ is $x$-even. Hence terms of $T$ and $D$ cannot cancel with each other, and both $T$ and $D$ vanish identically. Then one can proceed as in the proof of theorem \ref{main}. 

If $ii)$ holds, then $T$ is $y$-odd, $D$ is $y$-even.  Then one can proceed as above.
\hfill$\clubsuit$

\end{document}